\documentclass[12pt,a4paper]{article}
\usepackage[latin1]{inputenc}
\usepackage{graphicx}
\usepackage[mathscr]{euscript}
\usepackage{a4}
\usepackage{epsfig}
\usepackage{stmaryrd,amsfonts,amsmath,amssymb}

\newcommand{\origforall}{}    \let\origforall=\forall
\renewcommand{\forall}{\,\origforall\,}
\newcommand{\origexists}{}    \let\origexists=\exists
\renewcommand{\exists}{\,\origexists\,}
\newcommand{\origin}{}        \let\origin=\in
\renewcommand{\in}{{\,\origin\,}}
\newcommand{\origwedge}{}     \let\origwedge=\wedge
\renewcommand{\wedge}{{\quad\origwedge\quad}}
\newcommand{\origvee}{}       \let\origvee=\vee
\renewcommand{\vee}{{\quad\origvee\quad}}
\newcommand{\origfrac}{}      \let\origfrac=\frac
\renewcommand{\frac}[2]{{\;\origfrac{#1}{#2}\;}}
    
\renewcommand{\subset}{\subseteq}

\newcommand{\beq}{\begin{eqnarray*}}
\newcommand{\eeq}{\end{eqnarray*}}

\newcommand{\R}{\mathbb{R}}

\newcommand{\N}{\mathbb{N}}

\newcommand{\M}{\mathbb{M}}
\newcommand{\n}{\textnormal}

\DeclareMathOperator{\Lip}{{\n{Lip}}}
\newcommand{\smax}{\origvee}
\newcommand{\smin}{\origwedge}
\newcommand{\bigsup}{\bigvee}
\newcommand{\biginf}{\bigwedge}

\newcounter{claim}
\newtheorem{Theorem}[claim]{Theorem}
\newtheorem{Lemma}[claim]{Lemma}
\newtheorem{Corollary}[claim]{Corollary}
\newtheorem{Proposition}[claim]{Proposition}

\newtheorem{Example}[claim]{Example}
\newtheorem{Definition}[claim]{Definition}

\newcommand{\esquare}{\hfill\ensuremath{\square}}
\newenvironment{Proof}{{\textbf{Proof}$\;$}}{\esquare\newline}

\begin{document}

\title{Rough Isometries of Lipschitz Function Spaces}
\author{Andreas Lochmann}
\maketitle

\abstract

We show that rough isometries between metric spaces $X, Y$ can be lifted to
the spaces of real valued 1-Lipschitz functions over $X$ and $Y$ with supremum
metric and apply this to their scaling limits. For the inverse, we show how
rough isometries between $X$ and $Y$ can be reconstructed from structurally
enriched rough isometries between their Lipschitz function spaces.

\section{Introduction}

Is there a qualitative difference between functions on a continuum $X$ and
functions on a discrete set $Y$, which is $\epsilon$-dense in $X$? Obviously,
singularities may emerge on $X$. But, apart from that, are there more
differences? In this paper we want to show that, when we are concerned with
1-Lipschitz functions, nothing new is added in the passage from discrete to
continuum. Even more, the coarse geometry of a metric space (in the sense of
its rough isometries) is already determined by its 1-Lipschitz function space.

There are detailed books and lots of articles on many aspects of Lipschitz
function spaces, including isometries between them. A nice survey is Weaver's
book on {\it Lipschitz Algebras} \cite{Weaver} (cf. section 2.6, where Weaver
elaborates on the exact same questions we want to tackle here, but in a
non-coarse context and under different conditions). However, to the know\-ledge
of the author, no book nor paper dealt with their coarse geometry yet. On the
other hand, Lipschitz functions naturally appear in many aspects of coarse
geometry, like the Levy concentration phenomenon or the definition of
Lipschitz-Hausdorff distance in \cite{Gromov_Riemannian}. But they are not dealt
with as metric spaces either. 

The two main theorems we want to show are as following:

\begin{Theorem}\label{THEP_ind_ml-iso}
Let $X, Y$ be (possibly infinite) metric spaces, $\epsilon \geq 0$. For each
$\epsilon$-isometry $\eta: X\rightarrow Y$, there is a
$4\epsilon$-ml-isomorphism $\kappa: \Lip Y \rightarrow \Lip X$ such that
$\kappa$ is $\epsilon$-near $f\mapsto f\circ \eta$ for all $f\in \Lip Y$.
\end{Theorem}

\begin{Theorem}\label{THEP_ind_isom}
Let $X, Y$ be complete (possibly infinite) metric spaces and $\epsilon \geq 0$.
For each $\epsilon$-ml-isomorphism $\kappa: \Lip Y\rightarrow \Lip X$ there
is a $88\epsilon$-isometry $\eta: X\rightarrow Y$, such that $\kappa$
is $62\epsilon$-near $f\mapsto f\circ \eta$ for
all $f\in \Lip Y$. 
\end{Theorem}

The theorems above can be seen as variants on two 60 year old theorems. The
first is given by Hyers and Ulam in \cite{Hyers_Ulam}: 

\begin{Theorem}[Hyers-Ulam]\label{THE_hyers-ulam}
Let $K, K'$ by compact metric spaces and let $E$ and $E'$ be the spaces of all
real valued continuous functions on $K$ and $K'$, respectively. If $T(f)$ is a
homomorphism of $E$ onto $E'$ which is also an $\epsilon$-isometry, then there
exists an isometric transformation $U(f)$ of $E$ onto $E'$ such that
$||U(f)-T(f)|| \leq 21\epsilon$ for all $f$ in $E$. Corollary: The underlying
metric spaces $K$ and $K'$ are homeomorphic.
\end{Theorem}

As we make use of the order-lattice structure of the
Lipschitz function spaces in our proofs, the following theorem by Kaplansky 
\cite{Kaplansky} is related to our results as well, we state it 
in its formulation by Birkhoff (\cite{Birkhoff}, 2nd ed. p. 175f.):

\begin{Theorem}[Kaplansky]\label{THE_kaplansky}
Any compact Hausdorff space is (up to homeomorphism) determined by the
lattice of its continuous functions.
\end{Theorem}

The structure of this paper is as follows: In the first section we want to
recall the notions of rough isometry, lattices, and Lipschitz functions, and
our notation for metrics with infinity. In the second section we will proof
some simple properties of Lipschitz function spaces and define what we call
``$\Lambda$-functions''. In the third and fourth section, we will proof
the two main theorems of this paper. The fifth section presents an application
in the context of scaling limits, giving a partial answer to the introductory
question. The final section provides some conclusions and perspectives for
future development of this subject.

\subsection{Notation}

Throughout this paper, let $Z:=\R_{\geq 0}\cup \{\infty\}$ and $X, Y$ be
(possibly infinite) metric spaces in the following sense:

\begin{Definition}
A {\em (possibly infinite) metric space} $(X,d)$ is a non-empty set $X$
together with a mapping $d: X\times X \rightarrow Z$ which is
positive-definite and fulfills the triangle-inequality. For this, we set
$\infty + z = z + \infty = \infty$ and $z \leq \infty$ for all $z\in Z$ and
call $z\in Z$ positive iff $z\neq 0$. $(X,d)$ is called {\em true metric
space} iff $d(x,x')\neq\infty$ for all $x,x'\in X$.
\end{Definition}

We will make heavy use of the symbol
\beq
\big|d(x,x') - d(y,y')\big| &\leq& z
\eeq
for some $x,x',y,y'$ in $X$ or $Y$ and $z\in Z$. To make sense of this in case
one of the distances becomes infinite, we define the former symbol to be
equivalent to
\beq
d(x,x') \; \leq \; d(y, y') + z &\n{ and }& d(y, y') \; \leq \; d(x,x') + z.
\eeq
In particular, we find $|\infty - \infty| = 0$. This might seem
unfamiliar. Note however, that $|z - z'|$ can be perfectly understood as
a (possibly infinite) metric on $Z$ itself. (Note that there is no
non-trivial convergence to $\infty$ in this metric, $\infty$ is just
an infinitely far away point.)

In most cases we call metrics on $X$ and $Y$ both ``$d$'' as it should
be clear from the elements which metric is meant.

Furthermore, it's obvious that a (possibly infinite) metric space $X$
always is a disjoint union of true metric spaces $X_j$ with $d(x,x') =
\infty$ iff $x\in X_j, x'\in X_k$ with $j\neq k$, $j,k\in J$. We call
the $X_j$ {\em components} of $X$. We call $X$ {\em complete}, iff all
of its components are complete as true metric spaces.

\begin{Definition}
Two (set theoretic) mappings $\alpha,\beta: X\rightarrow Y$ are {\em
$\epsilon$-near} to each other, $\epsilon \geq 0$, iff $d_Y(\alpha x, \beta x)
\leq \epsilon \forall x\in X$. (We drop brackets where feasible.)

A (set theoretic) mapping $\alpha: X\rightarrow Y$ is {\em
$\epsilon$-surjective}, $\epsilon \geq 0$, iff for each $y\in Y$ there is
$x\in X$ s.t. $d_Y(\alpha x, y) \leq \epsilon$.
\end{Definition}

\begin{Definition}\label{DEF_coarse_isometry}
A (not neccessarily continuous) map $\eta: X\rightarrow Y$ is called an {\em
$\epsilon$-isometric embedding}, $\epsilon \geq 0$ (which shall always imply
$\epsilon\in\R$), iff
\beq
|d_X(x,x') - d_Y(\eta x, \eta x')| &\leq& \epsilon
\eeq
for all $x, x'\in X$. 

A pair $\eta: X\rightarrow Y$, $\eta': Y\rightarrow X$ of $\epsilon$-isometric
embeddings is called an {\em $\epsilon$-isometry} (or {\em rough isometry}) iff
$\eta\circ\eta'$ and $\eta'\circ \eta$ are $\epsilon$-near the identities on $Y$
and $X$, respectively. When we speak of an ``$\epsilon$-isometry
$\eta:X\rightarrow Y$'' a corresponding map $\eta'$ shall always be implied.

$X$ and $Y$ are called $\epsilon$-isometric, iff there is an
$\epsilon$-isometry between them. 
\end{Definition}

It's difficult to attribute the concept of rough isometry to a single person,
as it was always present in the notion of quasi-isometry, which itself was an
obvious generalization of what was then called pseudo-isometry by Mostow in his
1973-paper about rigidity (see \cite{Mostow}, \cite{Gromov_Hyperbolic},
\cite{Kanai}). Recent developments about the stability of rough isometries can
be found in \cite{Rassias_Isometries}.

\begin{Definition}
A mapping $f: X\rightarrow Y$ is called {\em $(K, \epsilon)$-Lipschitz}
(i.e. ``$K$-Lipschitz map on $\epsilon$-scale'' in \cite{Gromov_Asymptotic}),
$\epsilon, K\geq 0$, iff 
\beq
d_Y(f(x), f(y)) &\leq& K\cdot d_X(x,y) + \epsilon \qquad \forall
x,y\in X.
\eeq
If $\epsilon=0$, $f$ is {\em $K$-Lipschitz (continuous)}. Define
$\Lip_{K,\epsilon}(X,Y)$ to be the set of all $(K,\epsilon)$-Lipschitz
functions $X\rightarrow Y$, and $\Lip_{K,\epsilon} X :=
\Lip_{K,\epsilon}(X,Z)$, $\Lip X := \Lip_{1,0}(X)$.
\end{Definition}

If nothing else is said, $Z$ is the default target space for a Lipschitz
function.

Assume $f$ to be a $(K,\epsilon)$-Lipschitz function on $X$ and $f(x) =
\infty$ for some $x\in X$. Then clearly $f(y) = \infty$ for all $y$ in
finite distance to $x$. Thus, if $X$ is a true metric space, we have $\Lip X 
= \Lip(X,\R_{\geq 0}) \cup \{\infty\}$.

\begin{Definition}
A {\em lattice} $(L,\smin,\smax)$ is a set $L$ together with two mappings
$\smin, \smax: L\times L \rightarrow L$ which are commutative, associative and
fulfill the absorption laws $f\smin(f\smax g) = f\smax(f\smin g) = f$ for all
$f,g\in L$. A lattice $L$ is called {\em complete} iff all infima and all
suprema of all subsets of $L$ exist in $L$.
\end{Definition}

Of particular interest is $L=\Lip X$ with $\smin$ and $\smax$ pointwise
minimum and maximum respectively, and $\biginf$, $\bigsup$ pointwise
infimum and supremum. The following proposition is a special case of Lemma 6.3
in \cite{Heinonen_Analysis} and Proposition 1.5.5 in \cite{Weaver}. To keep
this article self-contained, we nevertheless give a proof:

\begin{Proposition}
Let $X$ be a (possibly infinite) metric space. Then $\Lip X$ is complete as a
lattice.
\end{Proposition}
\begin{Proof}
Let $f_j$, $j\in J$ be in $\Lip X$. Obviously, $Z$ is complete as a
lattice, with $\biginf_\emptyset = \infty$ and $\bigsup_\emptyset = 0$. So we
define pointwise 
\beq
g(x) \,:=\, \bigsup_{j\in J} f_j(x), &\quad & h(x) \,:=\, \biginf_{j\in J}
f_j(x)
\eeq
and observe that $g$ and $h: X\rightarrow Z$ are Lipschitz: Let $x,y\in X$ be
arbitrary. Then holds
\beq
h(x)\,\leq\,f_j(x)\,\leq\,f_j(y) + d(x,y)
\eeq
for all $j\in J$, and thus, by passing to the infimum:
\beq
h(x) &\leq& h(y) + d(x,y).
\eeq
Same for $g$.
\end{Proof}

\begin{Example}
The space $C([0,1], [0,1])$ of real-valued continuous functions
$[0,1]\rightarrow [0,1]$ is not a complete lattice with pointwise minimum and
maximum: Choose $f_j(x) = 0 \smax (1-j\cdot x)$, $j\in \N$, the infimum is not
continuous. The same example shows that the space $\bigcup_{K\geq 0}
\Lip_{K,0} X$ of all Lipschitz-functions with arbitrary Lipschitz constant
is no complete lattice.
\end{Example}

On $\Lip X$, we consider the (possibly infinite) supremum metric
\beq
d_\infty(f,g) &:=& \bigsup_{x\in X}\big|f(x)-g(x)\big|.
\eeq
Note that $(\Lip X, \smin, \smax, d_\infty)$ is no metric lattice in the sense
of Birkhoff \cite{Birkhoff}: There is no valuation on $\Lip X$ inducing
$d_\infty$, and property (4) of Theorem 1, p. 230 (third edition)
is explicitly violated even by bounded Lipschitz functions. 

\section{Fundamental properties}

\begin{Proposition}\label{PRO_ml}
For $f_j, g_j$ arbitrary set theoretic functions $X\rightarrow Z$, $j\in J$,
$J$ some arbitrary index set, holds:
\beq
d_\infty\left(\biginf_{j\in J} f_j, \biginf_{j\in J} g_j\right) &\leq&
\bigsup_{j\in J} d_\infty(f_j, g_j)\\
d_\infty\left(\bigsup_{j\in J} f_j, \bigsup_{j\in J} g_j\right) &\leq&
\bigsup_{j\in J} d_\infty(f_j, g_j)
\eeq
\end{Proposition}
\begin{Proof}
For $J=\emptyset$, both inequalities are trivial. Assume $J\neq \emptyset$.
As $\bigsup_{j\in J}$ and $\bigsup_{x\in X}$ commute, it suffices to show
\beq
d_\infty\left(\biginf x_j, \biginf y_j\right) &\leq&
\bigsup d_\infty(x_j, y_j)\\
d_\infty\left(\bigsup x_j, \bigsup y_j\right) &\leq& \bigsup
d_\infty(x_j, y_j)
\eeq
for any $x_j, y_j\in Z$.

First we handle infinities. First inequality: Assume there is $j$ with
$x_j=y_j=\infty$. We can ignore all such $j$'s from $J$, unless all
$x_j$ and $y_j$ are $\infty$. In this case on both sides are zeros.
Now assume $x_j=\infty\neq y_j$. Then $\infty$ appears on the right
side and trivializes the inequality. So we can restrict to finite
$x_j$ and $y_j$. Note that $\biginf_j x_j = \infty$ can only happen
when all $x_j = \infty$.

Second inequality: Assume $\bigsup_j x_j=\infty$, but $\bigsup_j y_j$
is finite. Then there is an upper bound for $y_j$ but not for
$x_j$. Hence the right side becomes infinite, too. Note that infinite
$x_j$ or $y_j$ automatically lead to infinite $\bigsup_j x_j$ or
$\bigsup_j y_j$, respectively.

Without restriction let $\biginf_j x_j \geq \biginf_j y_j$, and let $M
:= \bigsup_j d(x_j, y_j)$. Let $\delta > 0$ be arbitrary. Then there
is an $m\in J$ with $y_m \leq \biginf_j y_j + \delta$. Furthermore we
have $d(x_m, y_m)\leq M$, hence $y_m \geq x_m - M$. Altogether:
\beq
\biginf x_j \;\leq\; x_m \;\leq\; \biginf y_j + M + \delta
\eeq
Now let $\delta \rightarrow 0$.
The other inequality works the same way.
\end{Proof}

Next we define a special version of rough isometry, suiting the
lattice structure of Lipschitz function spaces. The main new property
will be an ``approximate lattice homomorphism''. It exists in various
versions, as Thomas Schick pointed out to us:

\begin{Proposition}\label{PRO_prop_ml_equiv}
Let $X, Y$ be (possibly infinite) metric spaces and $\kappa: \Lip Y \rightarrow
\Lip X$ an $\epsilon$-isometric embedding, $\epsilon \geq 0$. Then the following
properties are equivalent:
\begin{enumerate}
\item
$f\,\leq\,g \;\Rightarrow \;\kappa f \,\leq \, \kappa g + \epsilon$ for all
$f,g\in \Lip Y$
\item
$d_\infty\big((\kappa f) \smax (\kappa g), \kappa(f\smax g)\big) \,\leq\,
\epsilon$ for all $f,g\in \Lip Y$
\item For all $f_j\in\Lip Y$, $j\in J$, $J\neq \emptyset$ some index set, holds:
\beq
d_\infty \left(\bigsup_{j\in J} \kappa f_j,\; \kappa \bigsup_{j\in J}
f_j\right) \;\leq\; \epsilon &\n{and}& 
d_\infty \left(\biginf_{j\in J} \kappa f_j,\; \kappa \biginf_{j\in J}
f_j\right) \;\leq\; \epsilon
\eeq
\end{enumerate}
\end{Proposition}
\begin{Proof}
{\bf (3) $\Rightarrow$ (2):} $\# J = 2$.

{\bf (2) $\Rightarrow$ (1):} Assume there is some $x\in X$ such that $(\kappa
f)(x) \,>\, (\kappa g)(x) + \epsilon$. From $f\leq g$ follows $f\smax g = g$,
thus $d_\infty((\kappa f) \smax (\kappa g), \kappa g) \leq \epsilon$, in
particular $(\kappa f)(x) \smax (\kappa g)(x) \leq (\kappa g)(x) + \epsilon$,
contradiction.

{\bf (1) $\Rightarrow$ (3):} Obviously $f_k \leq \bigsup_{j\in J} f_j$ for all
$k\in J$, hence $\kappa f_k \leq \kappa \bigsup_{j} f_j + \epsilon$. We
calculate the supremum over all $k\in J$: $\bigsup_j \kappa f_j \leq
\kappa\bigsup_j f_j + \epsilon$. On the other hand, for every $\delta > 0$
there is some $k\in J$ with $d_\infty(f_k, \bigsup_j f_j) \leq \delta$. As
$\kappa$ is an $\epsilon$-isometric embedding, this yields $d_\infty(\kappa
f_k, \kappa \bigsup_j f_j) \leq \delta + \epsilon$, hence
\beq
\kappa \bigsup_{j\in J} f_j \;\leq\; \kappa f_k + \epsilon + \delta \;\leq\;
\bigsup_{j\in J} \kappa f_j + \epsilon + \delta.
\eeq
Let $\delta \rightarrow 0$. The other approximation works analogously.
\end{Proof}

We extend property (3) from Proposition \ref{PRO_prop_ml_equiv} to allow
$J=\emptyset$, and use it to define the notion of ml-isomorphisms:

\begin{Definition}\label{DEF_ml-iso}
Let $X, Y$ be (possibly infinite) metric spaces. An {\em
$\epsilon$-ml-homo\-morphism}, $\epsilon \geq 0$ is an $\epsilon$-isometric
embedding $\kappa: \Lip Y \rightarrow \Lip X$, with
\beq
d_\infty \left(\bigsup_{j\in J} \kappa f_j,\; \kappa \bigsup_{j\in J}
f_j\right) \;\leq\; \epsilon &\n{and}& 
d_\infty \left(\biginf_{j\in J} \kappa f_j,\; \kappa \biginf_{j\in J}
f_j\right) \;\leq\; \epsilon
\eeq
for all $f_j\in \Lip Y$, $j\in J$, $J$ some index set.

An {\em $\epsilon$-ml-isomorphism} is a pair of
$\epsilon$-ml-homomorphisms $\kappa: \Lip Y \rightarrow \Lip X$ and $\kappa':
\Lip X \rightarrow \Lip Y $, s.t. $\kappa\circ \kappa'$ and $\kappa'\circ
\kappa$ are $\epsilon$-near their corresponding identities. When we speak of an
``$\epsilon$-ml-isomorphism $\kappa$'', the corresponding $\kappa'$ shall
always be implied.
\end{Definition}

\begin{Proposition}\label{PRO_ml-iso_zero}
For any $\epsilon$-ml-isomorphism $\kappa$ holds $d_\infty(\kappa(0), 0)
\,\leq\, \epsilon$.
\end{Proposition}
\begin{Proof}
As Andreas Thom pointed out, this follows directly from Definition
\ref{DEF_ml-iso} when $J=\emptyset$. There's also a $3\epsilon$-proof
avoiding empty index sets:
Let $\kappa:\Lip Y\rightarrow \Lip X$ be an $\epsilon$-ml-isomorphism.
We certainly know $0\smin \kappa'(0) = 0\in \Lip Y$, hence
\beq
d_\infty(\kappa(0),\; \kappa(0)\smin \kappa\kappa'(0)) &\leq & 2\epsilon.
\eeq
Now apply Proposition \ref{PRO_ml} to see
\beq
d_\infty(\kappa(0)\smin\kappa\kappa'(0),\; \kappa(0)\smin 0) &\leq& \epsilon
\eeq
and use $\kappa(0)\smin 0 = 0\in\Lip X$.
\end{Proof}

\begin{Proposition}
A $\delta$-surjective $\epsilon$-ml-homomorphism $\kappa: \Lip Y\rightarrow
\Lip X$ induces a $(2\epsilon+2\delta)$-ml-isomorphism $(\kappa, \kappa')$.
\end{Proposition}
\begin{Proof}
For each $f\in \Lip X$ choose an element $\kappa'(f)$,
s.t. $d(\kappa\kappa'f, f) \leq \delta$. 
We show that the pair $(\kappa,  \kappa')$ defines a
$(2\epsilon+3\delta)$-ml-isomorphism. The first inequality in Definition
\ref{DEF_ml-iso} is standard in coarse geometry:
\beq
d_\infty(f,g) \;\leq\; d_\infty(\kappa\kappa'f,\kappa\kappa'g) +
2\delta \;\leq\; d_\infty(\kappa'f, \kappa'g) + (\epsilon + 2\delta)\\
d_\infty(f,g) \;\geq\; d_\infty(\kappa\kappa'f,\kappa\kappa'g) -
2\delta \;\geq\; d_\infty(\kappa'f, \kappa'g) - (\epsilon + 2\delta)
\eeq
for all $f,g\in \Lip X$. We now show that $\kappa'$ fulfills the second and
third inequality as well. Both can be handled the same way:
\beq
d_\infty \left(\biginf \kappa' f_j,\; \kappa' \biginf f_j\right) &\leq&
d_\infty\left(\kappa \biginf \kappa' f_j, \; \kappa \kappa' \biginf f_j\right) +
\epsilon \\ &\leq& d_\infty \left(\kappa \biginf \kappa' f_j, \; \biginf f_j
\right) + \epsilon + \delta \\ &\leq& d_\infty \left(\biginf \kappa \kappa' f_j,
\; \biginf f_j \right) + 2\epsilon + \delta\\ &\leq& \bigsup d_\infty
\left(\kappa \kappa' f_j, \; f_j\right) + 2\epsilon + \delta\\ &\leq& \delta +
2\epsilon + \delta \qquad \forall f_j\in \Lip X, j\in J
\eeq
Here we used (i) $\kappa$ is $\epsilon$-isometric embedding, (ii)
$\kappa\kappa'$ is near identity, (iii) $\kappa$ is ml-homomorphism, (iv)
Proposition \ref{PRO_ml}, (v) $\kappa\kappa'$ is near identity.

Finally we show that $\kappa'\kappa$ is $(\epsilon+\delta)$-near identity:
\beq
d_\infty(\kappa'\kappa f, f) \;\leq\; d_\infty(\kappa\kappa'(\kappa
f), (\kappa f)) + \epsilon \;\leq\; \delta + \epsilon \qquad \forall
f\in \Lip X.
\eeq
\end{Proof}

$\Lip X$ is no algebra, like e.g. $C(X)$. Thus we can't give a basis of
functions and reconstruct $\Lip X$ by linear combinations. However, we can use
the lattice structure to give another kind of ``basis'' for $\Lip X$: Minimal
Lipschitz functions with a given value at a single point.

\begin{Definition}\label{DEF_lambda_function}
Let $x, y\in X$ and $r\in Z$ be arbitrary. Define $\Lambda(x,r)\in
\Lip X$ by $\Lambda(x,r)(y) := (r - d(x,y)) \smax 0$.
\end{Definition}

Note that this definition applies to $r=\infty$ or $d(x,y)=\infty$ as well: If
$d(x,y)=\infty$, we have $\Lambda(x,r)(y)=0$, and if $r=\infty$:
\beq
\Lambda(x,\infty)(y) &=&\left\{\begin{array}{lcl} \infty &:& d(x,y) \neq
    \infty \\ 0 &:& d(x,y) = \infty\end{array}\right.
\eeq
$\Lambda$-functions with $r=\infty$ will be called {\em infinite}, else {\em
finite}. Infinite $\Lambda$-functions are infinitely high characteristic
functions for $X$'s components. 

\begin{Proposition}\label{PRO_lambda_distances}
Let $x,y\in X$, $r,s\in Z$. Then holds:
\beq
d_\infty(\Lambda(x,r), \Lambda(y,s)) &=& 
\left\{\begin{array}{lcl}
r\smax s &:& d(x,y) \,\geq\, r\smin s\\
|r-s|+d(x,y) &:& d(x,y) \,\leq\, r\smin s \,<\, \infty\\
0&:& d(x,y) \,<\, r\smin s \,=\, \infty
\end{array}\right.\\
&\leq& |r-s| + d(x,y)
\eeq
\end{Proposition}
\begin{Proof}
Note that if $d(x,y)=r\smin s$ the first and second case coincide, as $|r-s| =
r\smax s - r\smin s$. Assume without restriction $r\leq s$. Let
\beq
fz &:=& \big|\big(0\smax(r-d(x,z))\big) - \big(0\smax (s-d(y,z))\big)\big|
\quad \forall z\in X\\
d &:=& d_\infty(\Lambda(x,r),\Lambda(y,s)) \;=\; \bigsup_{z\in X} f(z).
\eeq

Let's start with infinite cases. If $r = s = d(x,y) = \infty$, we get
$d=\infty$ on both sides. If $r = s =\infty$, $d(x,y) \neq \infty$, we get
$d=0$. This is correct, as in this case the $\Lambda$-functions are equal. If 
$s=\infty$, $r\neq \infty$ we get $d=\infty$ again, for each variant of
$d(x,y)$. If $r,s\neq \infty$ but $d(x,y)=\infty$, the two $\Lambda$-functions
have different components as support, and thus $d$ becomes the maximum of the
differences, this is $s$.

Now we assume $r,s,d(x,y)\neq \infty$. First case: $r\leq d(x,y)$. Then we
have
\beq
d &\geq& fy \;=\; \big|s - \big(0\smax(r-d(x,y))\big)\big| \; = \; s
\eeq
In addition, we have $\Lambda(x,r)(z) \in [0,r]$, $\Lambda(y,s)(z)\in [0,s]$,
thus $fz\leq r\smax s = s$. Hence $d=s$. Second case: $d(x,y)\leq r$.
\beq
d &\geq& fy \;=\; \big|s - \big(0\smax(r-d(x,y))\big)\big| \; = \;
s-r+d(x,y)
\eeq
And:
\beq
fz &=& |r-d(x,z)-s+d(y,z)| \;\leq\; |r-s| + |d(y,z) - d(x,z)|\\
&\leq& |r-s| + d(x,y) \;=\; s - r + d(x,y) \quad \forall z\in X
\eeq
\end{Proof}

\begin{figure}[t]
\centerline{\epsfig{file=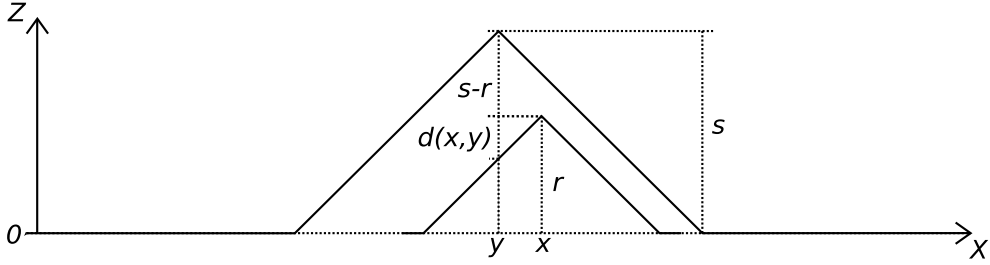}}
\caption{The $d_\infty$-distance between two $\Lambda$-functions is determined
by their difference evaluated at the maximum point of the larger function, see
Prop. \ref{PRO_lambda_distances}.}
\label{FIG_lambda_distance}
\end{figure}

\begin{Corollary}\label{COR_lambda_metric}
For all $x, y\in X$ holds:
\beq
d(x,y) &=& \lim_{r\rightarrow \infty,\, r\neq \infty}
d_\infty\big(\Lambda(x,r),\Lambda(y,r)\big)
\eeq
\end{Corollary}
\begin{Proof}
Follows directly from Proposition \ref{PRO_lambda_distances}.
\end{Proof}

This Corollary points us at an interesting aspect of $\Lambda$-functions: When
we analyse the metric space $X_r := \{\Lambda(x,r):x\in X\}$ with metric
$d_\infty$ for a fixed $r \in \R_{>0}$, we find it naturally isometric
to $(X, d_r)$ with the cut-off-metric $d_r(x,y) := r \smin d(x,y)$ for
all $x,y\in X$. Only in the limit $r\rightarrow \infty$, $d_\infty$
will restore the full metric of $X$. Ironically, $d_\infty$ obviously
cuts away the coarse, large-scale information of $X$ (in which we're
primarily interested) and conserves the topological, small-scale
informations. The large-scale information of $X$ is still present, but
more subtle to access.

\begin{Proposition}\label{PRO_lambda_building}
For all $f\in\Lip X$ holds: $f = \bigsup_{x\in X} \Lambda(x, f(x))$, where the
latter is a pointwise maximum, not only supremum. If $X$ is complete, the set
$A := \{\Lambda(x,f(x)):x\in X\}\cup\{0\}$ is (topologically) closed.
\end{Proposition}
\begin{Proof}
Let $g := \bigsup_{x\in X}\Lambda(x,fx)$. Clearly, we have $\forall z\in X:$
$fz \leq gz$, as $fz=\Lambda(z,fz)(z)$. We now observe that
\beq
fz &\geq& \Lambda(x,fx)(z) \qquad \forall x, z\in X.
\eeq
For $d(x,z) \geq fx$, this is clear. For $d(x,z) \leq fx$, this follows
from Lipschitz continuity ($fz \;\geq\; fx - d(x,z)$).

Furthermore, we notice that we deal with pointwise maxima: Each supremum of
$\{\Lambda(x,fx)(z)\}_{x\in X}$ is taken by $\Lambda(z,fz)(z)=fz$.

Let $\Lambda(x_j, fx_j)$ be any sequence converging to $g\in\Lip X$, $x_j\in
X$, $j\in\N$. First we notice 
\beq
fx_j \;=\; d_\infty(0, \Lambda(x_j,fx_j)) &\geq& d_\infty(0, g) -
d_\infty(g,\Lambda(x_j, fx_j))\\
\n{and } fx_j &\leq& d_\infty(0,g) + d_\infty(g,\Lambda(x_j, fx_j))
\eeq
hence $fx_j \rightarrow d_\infty(0,g)$. Assume $g\neq 0$ and
finite. Then there is $x\in X$ with $gx > 0$ and $fx_j$ must have a
lower bound $R>0$ for large enough $j$. By Cauchy criterion there is
$N\in\N$ such that for all $j,k>N$ we have
\beq
d_\infty(\Lambda(x_j,fx_j), \Lambda(x_k,fx_k)) \;\leq\; \frac{1}{2}R \;<\;
fx_j\smin fx_k.
\eeq
Due to Proposition \ref{PRO_lambda_distances} we conclude that for
large enough $j,k$
\beq
d_\infty(\Lambda(x_j,fx_j), \Lambda(x_k,fx_k)) \;=\; |fx_j - fx_k| +
d(x_j, x_k) \rightarrow 0.
\eeq
Thus $fx_j$ as well as $x_j$ are Cauchy-sequences. As $X$ and
$Z$ are metrically complete, we find $x' := \lim x_j$. As $f$ is
continuous, we have $fx' = \lim fx_j$, and $\Lambda(x',fx')\in
A$. Now we only have to show $g = \Lambda(x', fx')$. But this is
clear, as for large enough $j$ we have
\beq
d_\infty(\Lambda(x',fx'), \Lambda(x_j, fx_j)) &\leq& |fx' - fx_j| +
d(x', x_j) \;\rightarrow\; 0+0.
\eeq
Now assume $g$ to be infinite (i.e. $\exists x: gx=\infty$). Then
$fx_j$ has to be infinite as well for large enough $j$ (there is no
non-trivial convergence to infinity in the chosen metric on $Z$) and
Proposition \ref{PRO_lambda_distances} shows $d(x_j, x_k)<\infty$ for
large enough $j,k$. Hence $\Lambda(x_j, fx_j) = \Lambda(x_k, fx_k) =
g$.
\end{Proof}

We make some more use of the black magic of Proposition \ref{PRO_ml}:

\begin{Proposition}\label{PRO_lambda_exchange}
For all $\epsilon$-isometries $\eta: X\rightarrow Y$ and $f\in \Lip Y$ holds:
\beq
d_\infty\left(\bigsup_{x\in X} \Lambda(x,f\eta x), \bigsup_{y\in Y}
\Lambda(\eta' y, fy)\right) &\leq& \epsilon
\eeq
\end{Proposition}
\begin{Proof}
We observe that $d := d_\infty(\bigsup_{x\in X} \Lambda(x,f\eta x), \bigsup_{y\in Y}
\Lambda(\eta' y, fy))$ can be rewritten to
\beq
d &=& d_\infty\left(\bigsup_{(x,y)\in J} \Lambda(x,f\eta x), \bigsup_{(x,y)\in J}
\Lambda(\eta' y, fy)\right)
\eeq
where $J := \{(x,y)\in X\times Y: y=\eta x \n{ or } x = \eta'y\}$:
Each element of $X$ (respectively $Y$) appears at least once in $J$, and
multiple instances don't matter, as $\bigsup$ is idempotent. Now Proposition
\ref{PRO_ml} yields:
\beq
d &\leq& \bigsup_{(x,y)\in J} d_\infty\big(\Lambda(x,f\eta x), \Lambda(\eta' y,
fy)\big)
\eeq
Let $(x,y)\in J$. Case 1: $y=\eta x$. Then 
\beq
d_\infty\big(\Lambda(x,f\eta x), \Lambda(\eta' y, fy)\big) &=&
d_\infty\big(\Lambda(x,f\eta x), \Lambda(\eta'\eta x, f\eta x)\big)\\
&\leq& d(x,\eta'\eta x) \;\leq\; \epsilon
\eeq
Case 2: $x=\eta' y$:
\beq
d_\infty\big(\Lambda(x,f\eta x), \Lambda(\eta' y, fy)\big) &=&
d_\infty\big(\Lambda(\eta'y,f\eta \eta'y), \Lambda(\eta'y, fy)\big)\\
&\leq& |f\eta\eta' y - fy| \;\leq\; \epsilon
\eeq
\end{Proof}

\section{Inducing rough ml-isomorphisms}

We make a first use of the notions of the preceding section. We proof that
each $\epsilon$-isometry $\eta:X\rightarrow Y$ lifts to an $\epsilon$-isometry
$\bar{\eta}: \Lip Y \rightarrow \Lip X$. Even better, $\bar{\eta}$ is an
$\epsilon$-ml-isomorphism, and is near $f\mapsto f\circ \eta$.

The next lemma is kind of a smoothening theorem. It states that the space of
$1$-Lipschitz functions over $X$ is $\epsilon$-dense in the space
of $(1,\epsilon)$-Lipschitz functions over $X$ for all $\epsilon \geq 0$. A
similar result for continuous functions is given by Petersen in
\cite{Petersen}, section 4.

\begin{Lemma}
\label{LEM_lipschitzisation}
Let $f\in\Lip_{1,\epsilon}(X,Z)$. Define
\beq
\bar{f} := \bigsup_{x\in X} \Lambda(x, fx)
\eeq
Then $f$ and $\bar{f}$ are $\epsilon$-near.
\end{Lemma}
\begin{Proof}
We observe that $f(y)$ is never larger than $\bigsup_{x\in
X}\Lambda(x,fx)(y)$ for all $y\in X$. So we have
\beq
d_\infty\left(f,\bigsup_{x\in X} \Lambda(x, fx)\right)&=&\bigsup_{x,y\in X}
\big(\Lambda(x,fx)(y)-f(y)\big)
\eeq
and furthermore
\beq
\Lambda(x,fx)(y) - f(y) &=& \left\{\begin{array}{ll} -f(y)&:
d(x,y)\,\geq\, f(x)\\ f(x)-f(y)-d(x,y) &: d(x,y) \,\leq\, f(x)
\end{array}\right..
\eeq
As $f(x)-f(y)-d(x,y) \leq \epsilon$ and $-f(y)\leq 0 \leq\epsilon$ we
conclude the statement. (Note that each negative value is surpassed by at
least one non-negative value, i.e. $-f(y)$ never occurs after taking the
supremum.)
\end{Proof}

\begin{Proposition}\label{PRO_lip2}
If $\eta: X\rightarrow Y$ is an $\epsilon$-isometry, and $\kappa:
\Lip Y \rightarrow \Lip X$ any mapping which is $\delta$-near $f\mapsto
f\circ \eta$, then $\kappa$ is a $(2\epsilon + 2\delta)$-ml-isomorphism.
\end{Proposition}
\begin{Proof}
(i) We show $|d_\infty(f, g) - d_\infty(\kappa f, \kappa g)| \leq 2\epsilon +
2\delta$ for all $f, g\in \Lip X$. We have
\beq
\big|d_\infty(\kappa f, \kappa g) - d_\infty(f\circ \eta, g\circ \eta)\big|
&\leq& 2\delta.
\eeq
As next we notice $d_\infty(f\circ\eta, g\circ\eta) \leq
d_\infty(f,g)$. Now let $x\in X$ be arbitrary, $y\in Y$ such that $d_Y(y,\eta
x)\leq \epsilon$. Then $|f\eta x - fy| \leq \epsilon$ as $f$ is
1-Lipschitz. Hence
\beq
&&|f\eta x - g\eta x| \;\leq\; |fy - gy| + 2\epsilon \;\leq\;
d_\infty(f, g) + 2\epsilon\\
\Rightarrow \quad && d_\infty(f\circ \eta, g\circ \eta) \;\leq\; d_\infty(f,g)
+ 2\epsilon.
\eeq

(ii) For $J=\emptyset$ we observe that $d_\infty(\kappa(0), 0\circ \eta)
\,\leq\, \delta$ and $0\circ \eta = 0$, as well as $d_\infty(\kappa(\infty),
\infty \circ \eta) \,\leq\, \delta$ and $\infty\circ\eta = \infty$. Hence,
assume $J\neq \emptyset$. We know
\beq
d_\infty\left(\kappa\left(\biginf f_j\right), \biginf (f_j\circ \eta)\right)
&=& d_\infty\left(\kappa\left(\biginf f_j\right), \left(\biginf f_j\right)
\circ \eta\right) \;\leq\; \delta
\eeq
as the infimum is calculated pointwise. Hence, with Proposition
\ref{PRO_ml}:
\beq
d_\infty\left(\biginf(\kappa f_j), \kappa\left(\biginf
f_j\right)\right) &\leq & d_\infty\left(\biginf(\kappa f_j), \biginf
(f_j\circ\eta)\right) + \delta\\ &\leq& \bigsup d_\infty(\kappa f_j,
f_j\circ \eta) + \delta\\ &\leq& \epsilon + \delta.
\eeq
Same for supremum.
\end{Proof}

\begin{Theorem}[$=$ Th. \ref{THEP_ind_ml-iso}]\label{THE_ind_ml-iso}
Given an $\epsilon$-isometry $\eta: X\rightarrow Y$, $\bar{\eta}(f) :=
\overline{f\circ \eta}$ defines a $4\epsilon$-ml-isomorphism from $\Lip Y$ to
$\Lip X$.
\end{Theorem}
\begin{Proof}
Let $f\in\Lip Y$ be arbitrary. $f\circ\eta$ satisfies
\beq
d(f\eta x, f\eta y) \;\leq\; d(\eta x, \eta y) \;\leq\; d(x,y) + \epsilon.
\eeq
Hence, $f\circ\eta$ and $\bar{\eta}(f)$ are $\epsilon$-near (Lemma
\ref{LEM_lipschitzisation}). However, $\bar{\eta}(f)$ is in $\Lip X$, as it is
a supremum of Lipschitz functions. Thus we can apply Proposition
\ref{PRO_lip2} to $\bar{\eta}:\Lip Y \rightarrow \Lip X$. Same holds for $\eta'$
(Definition \ref{DEF_coarse_isometry}). It remains to show that
$\bar{\eta}\circ \bar{\eta'}$ and $\bar{\eta'}\circ\bar{\eta}$ are near their
respective identities. 

We already saw that $\bar{\eta}(\bar{\eta'}(f))$ is $\epsilon$-near
$(\bar{\eta'}f)\circ \eta$. Similarly $\bar{\eta'}f$ is $\epsilon$-near
$f\circ \eta'$ and thus $(\bar{\eta'}f)\circ\eta$ is $\epsilon$-near
$f\circ\eta'\circ\eta$. Finally, $\eta'\circ\eta$ is $\epsilon$-near identity,
and as $f$ is 1-Lipschitz, $f\circ\eta'\circ\eta$ is $\epsilon$-near $f$,
too. All this adds up to $3\epsilon$. Same for $\bar{\eta'}\circ\bar{\eta}$.
\end{Proof}

\section{Inducing rough isometries}

In this section, we show the reversal of Theorem \ref{THE_ind_ml-iso}: Given
an $\epsilon$-ml-isomorphism $\kappa$ we construct a rough isometry $\eta$
such that $\bar{\eta}$ is near $\kappa$.

Recall the definition of a join-irreducible:
\beq
f = g \smax h &\Rightarrow& f = g \vee f = h
\eeq
It is interesting to see that the finite $\Lambda$-functions defined in
Definition \ref{DEF_lambda_function} satisfy a much more powerful version of
join-irreducibility:

\begin{Lemma}\label{LEM_lambda_property}
Let $p\in \Lip Y$, $Y$ complete. The following are equivalent:
\begin{enumerate}
\item $p$ is a finite $\Lambda$-function, i.e. $\exists x\in Y, r\in
Z\setminus\{\infty\}: p = \Lambda(x,r)$,
\item $\forall (f_j)_{j\in J}\subset \Lip Y, R \in Z:$
\beq
  d_\infty\left(p,\, \bigsup_{j\in J} f_j\right) \;\leq\; R
  &\Rightarrow& \forall \delta > 0 \exists j\in J: \; d_\infty (p, f_j) \;\leq\;
  R + \delta
\eeq
\end{enumerate}
\end{Lemma}
\begin{Proof}
In (2), the case $R=\infty$ is trivial. Hence, assume $R$ to be finite.

(1)$\Rightarrow$(2): Let $(f_j)_{j\in J}\subset \Lip Y$ and $R\geq 0$ be
   s.t. $d(p,\bigsup_{j\in J} f_j) \leq R$ holds. Choose $\delta > 0$
   arbitrary and $p=\Lambda(y,s)$, $y\in Y$, $s\in Z\setminus
   \{\infty\}$. As
\beq
   d\left(p(y), \bigsup f_j(y)\right) \;\leq\; R &\Rightarrow& p(y) - R -
   \delta \;<\; \bigsup f_j(y),
\eeq
   there has to be a $k\in J$ such that $p(y)-R-\delta < f_k(y)$, otherwise
   $p(y) - R-\delta$ would be a smaller upper bound for all $f_j$ then
   $\bigsup f_j(y)$. From this, we see
\beq
   f_k(x) &\geq& f_k(y) - d(x,y) \;>\; p(y) - d(x,y) - R - \delta.
\eeq
Case 1: $p(y) \geq d(x,y)$. Then we have $p(x) = p(y) - d(x,y)$, and
\beq
f_k(x) &>& p(x) - R - \delta.
\eeq
Case 2: $p(y) \leq d(x,y)$. Then $p(x)=0$ and
\beq
f_k(x) \;\geq\; 0 \;>\; p(x) - R - \delta
\eeq
holds trivially.

On the other hand, we have
\beq
f_k(x) &\leq& \bigsup_{j\in J} f_j(x) \;\leq\; p(x) + R \;<\; p(x) + R +
   \delta \qquad \forall x\in Y
\eeq
   and thus $d_\infty(f_k, p) < R + \delta$.

(2)$\Rightarrow$(1): Choose $J=Y$, $f_y=\Lambda(y,p(y))$, $R=0$, $\delta =
   1/n$. This yields a sequence $y_n$ of indizes (= points in $Y$) such that
   $\Lambda(y_n,p(y_n)) \rightarrow p$. As $\{\Lambda(y,p(y)):y\in
   Y\}\cup\{0\}$ is closed, we have either $p = \Lambda(y,p(y))$ for some
   $y\in Y$, or $p = 0 = \Lambda(y,0)$ for any $y\in Y$. Now assume
   $p(y)=\infty$. Then
\beq
   d_\infty\left(p, \bigsup_{r\in Z\setminus\{\infty\}}\Lambda(y,r)\right) \;
   &=& \; 0\\ \Rightarrow \exists r\in Z\setminus\{\infty\}: \;
   d_\infty\big(\Lambda(y,\infty), \Lambda(y,r)\big) &\leq& 1
\eeq
   This is a contradiction to Proposition \ref{PRO_lambda_distances}, hence $p$
   is a finite $\Lambda$-function.
\end{Proof}

\begin{figure}[t]
\centerline{\epsfig{file=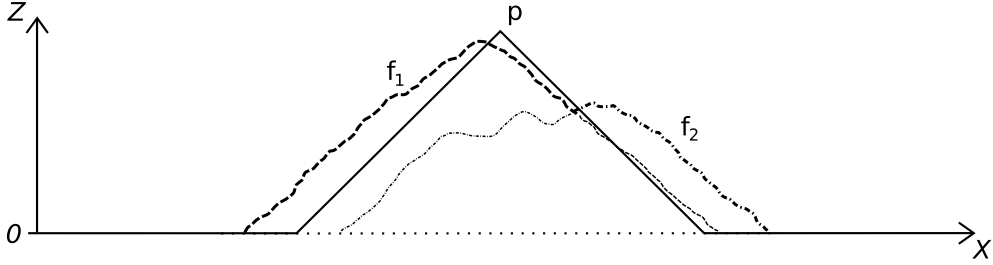}}
\caption{When approximating a $\Lambda$-function $p$ by Lipschitz functions
$f_j$, one of the functions (here $f_1$) must approximate the maximum point of
$p$. This function may not decrease too fast (Lipschitz!), and may not increase
too fast, as it is bounded from above by the approximation of $p$, hence it
already approximates $p$ on its own, see Lemma \ref{LEM_lambda_property}.}
\label{FIG_lambda_ml-property}
\end{figure}

\begin{Example}
Lemma \ref{LEM_lambda_property} does not hold for $\delta = 0$, just
insert $\Lambda(y, 1) = \bigsup_{r\in(0,1)} \Lambda(y, r)$.
\end{Example}

Recalling the short note after Corollary \ref{COR_lambda_metric}, the metric
information of $Y$ is encoded in the $\Lambda$-functions and the distances
between them. However, these functions are at first sight just some arbitrary
subset of $\Lip Y$ and thus there's no hope for the metric space
$(\Lip Y, d_\infty)$ to hold the full information about $Y$'s metric. The
preceding Lemma now explains to us, that the (finite) $\Lambda$-functions are
not arbitrary at all -- they have a specific, lattice theoretic property that
distinguishes them from the remaining functions. Hence, in some sense the metric
information of $Y$ is now part of the combined metric and lattice structure of
$\Lip Y$.

\begin{Proposition}\label{PRO_lambda_complete}
Let $Y$ be complete. Then the set $\{\Lambda(y,r):y\in Y, r\in Z\}$ of all
$\Lambda$-functions in $\Lip Y$ is topologically closed. In addition, the set
of all finite $\Lambda$-functions is topologically closed.
\end{Proposition}
\begin{Proof}
Let $(f_j)_{j\in\N} = (\Lambda(x_j, r_j))_{j\in \N}$ be some sequence of
$\Lambda$-elements in $\Lip Y$ with limit $f$. If there is a subsequence $j(n)$
with $r_{j(n)}\rightarrow 0$ for $n\rightarrow \infty$, then this subsequence
and hence $(f_j)$ converges to $f=\Lambda(x,0)$ for any $x\in Y$. So assume
$\inf r_j$ is positive for large enough $j$. Due to Proposition
\ref{PRO_lambda_distances} we have:
\beq
d(f_j, f_k) &=& \left\{\begin{array}{lcl}r_j \smax r_k&:&d(x_j,x_k) \geq
    r_j\smin r_k \\
|r_j-r_k|+d(x_j,x_k)&:&d(x_j,x_k) \leq r_j\smin r_k <\infty \\
0&:& d(x_j, x_k) < r_j \smin r_k = \infty
\end{array}
\right.
\eeq
As the left side becomes arbitrarily small, whereas
$r_j\smax r_k$ has a positive lower limit, only the second and third
case may occur for $j,k\rightarrow \infty$. For large enough $j,k$,
these cases don't mix anymore. The third case is trivial. From the
second case we conclude $|r_j - r_k|\rightarrow 0$ and $d(x_j, x_k)
\rightarrow 0$, and thus $r_j\rightarrow : r$ and $x_j\rightarrow :
x$. Clearly, $f=\Lambda(x,r)$. In particular, $r$ is finite in this
case, which proofs the second statement of the Proposition.
\end{Proof}

\begin{Lemma}\label{LEM_lambda_to_lambda}
Let $X$, $Y$ be complete, $\epsilon\geq 0$ and $\kappa: \Lip Y\rightarrow \Lip
X$ an $\epsilon$-ml-isomorphism. Then $\kappa$ maps finite $\Lambda$-functions
$6\epsilon$-near finite $\Lambda$-functions.
\end{Lemma}
\begin{Proof}
Let $p$ be some finite $\Lambda$-function. Represent $\kappa(p)$
via $\Lambda$-functions $q_j$, $j\in J$. Let $\delta > 0$ be arbitrary. Then we
have
\beq
d_\infty \left(\kappa (p), \bigsup q_j\right) &=& 0\qquad \qquad |\n{ apply
}\kappa' \\ \Rightarrow \qquad d_\infty \left(p, \bigsup \kappa'(q_j)\right)
&\leq& 3\epsilon
\eeq
Applying Lemma \ref{LEM_lambda_property} to $p$, we know that there exists
$k\in J$ such that
\beq
d_\infty (p,\, \kappa'(q_k)) &\leq& 3\epsilon + \delta\qquad \qquad |\n{ apply
}\kappa\\ \Rightarrow \qquad d_\infty (\kappa (p),\, q_k) &\leq& 5\epsilon +
\delta.
\eeq
$q_k$ must be a finite $\Lambda$-function, as
\beq
d_\infty(0, q_k) \;\leq\; d_\infty(0,\kappa(p)) + 5\epsilon + \delta \;\leq\;
d_\infty(0,p) + 8\epsilon + \delta \;<\;\infty.
\eeq

Case 1: $\epsilon > 0$. Choose $\delta = \epsilon$.

Case 2: $\epsilon = 0$. The preceding argument yields a sequence of finite
$\Lambda$-functions converging to $\kappa(p)$. As of Proposition
\ref{PRO_lambda_complete}, $\kappa(p)$ must be a finite $\Lambda$-function as
well.
\end{Proof}

The preceding Lemma is the critical point in our analysis: We can use
$\Lambda$-functions as building blocks for Lipschitz functions, as Proposition
\ref{PRO_lambda_building} tells us. From Lemma \ref{LEM_lambda_to_lambda} we
now know that these building blocks (or, at least, the finite versions) behave
sensible under $\epsilon$-ml-isomorphisms $\kappa$, such that we only have to
understand how they are mapped by $\kappa$ to reconstruct all other Lipschitz
functions. In particular, as they are strongly connected to the underlying
spaces, they allow us to define mappings between them:

\begin{Lemma}\label{LEM_mt1}
Let $X, Y$ be complete, $\epsilon \geq 0$.
Let $\kappa: \Lip Y \rightarrow \Lip X$ be an $\epsilon$-ml-isomorphism,
$\epsilon \geq 0$. Then there is a map $\eta: X\rightarrow Y$ such that
\beq
d_\infty\big(\Lambda(\eta x, r), \kappa'(\Lambda(x,r))\big) &\leq & 
59\epsilon
\eeq
for all $x\in X$, $r\in Z$. For $r\in [38\epsilon,\infty)$, we may replace
``$59\epsilon$'' by ``$43\epsilon$''.
\end{Lemma}
\begin{Proof}
In the following proof, the first two cases will deal with $\epsilon >
0$ and finite $r$, the third with $\epsilon=0$ and finite $r$ and the
fourth with $r=\infty$. 

{\bf Case 1 and 2:} For each $x\in X$, choose $\eta(x)\in Y$ and
$s_x\in Z\setminus \{\infty\}$ such that $\Lambda(\eta x, s_x)$ is
$6\epsilon$-near $\kappa'\Lambda(x, 22\epsilon)$ (use Lemma
\ref{LEM_lambda_to_lambda}).

{\bf Case 1:} $\epsilon > 0$, $r\in [38\epsilon, \infty)$. Let
$\Lambda(x',r')$ be $6\epsilon$-near $\kappa'\Lambda(x,r)$. Then by
Proposition \ref{PRO_ml-iso_zero} holds
\beq
d_\infty(0,\; \Lambda(x,r)) \;=\; r \quad\Rightarrow\quad
\big|d_\infty(0,\;\kappa'\Lambda(x,r)) - r\big| &\leq& 2\epsilon\\
\Rightarrow\qquad |r' - r| &\leq& 8\epsilon.
\eeq
In the same way, we have
\beq
\big|d_\infty(0,\; \Lambda(\eta x, s_x)) -
d_\infty(0,\; \kappa'\Lambda(x,22\epsilon))\big| &\leq& 6\epsilon\\
\Rightarrow \qquad \big|s_x - 22\epsilon\big| &\leq& 8\epsilon.
\eeq
We now take a look at
\beq
d_\infty(\Lambda(x,r),\;\Lambda(x,22\epsilon)) &=& r - 22\epsilon \quad
\n{(as $r\geq 22\epsilon$)}\\
\Rightarrow \qquad \big|d_\infty(\kappa'\Lambda(x,r),\; \kappa'\Lambda(x,
22\epsilon)) - (r - 22\epsilon)\big| &\leq& \epsilon\\
\Rightarrow \qquad \big|d_\infty(\Lambda(x',r'),\;\Lambda(\eta x, s_x)) -
(r-22\epsilon)\big| &\leq& 13\epsilon.
\eeq
Now we calculate $d := d_\infty(\Lambda(x',r'),\Lambda(\eta x, s_x))$ by hand.
 From Proposition \ref{PRO_lambda_distances}, $d$ could be $r'\smax s_x$ or
$d(x',\eta x) + |r'-s_x|$. We know
\beq
s_x \;\leq\; 8\epsilon + 22\epsilon \;=\;
30\epsilon \;\leq\; r - 8\epsilon 
\;\leq\; r',
\eeq
hence $r'\smax s_x = r'$. But, as $d \leq r - 22\epsilon + 13\epsilon
= r - 9\epsilon$, but $r' \geq r - 8\epsilon$, $d$
can't be $r'$ (here we use $\epsilon > 0$).
Remains 
\beq
d \;=\; d(x',\eta x) + |r'-s_x| \qquad \n{with}\qquad |d - (r-22\epsilon)|
\;\leq\; 13\epsilon.
\eeq
As shown above, $r' \geq s_x$, hence
\beq
d(x',\,\eta x) &\leq& r - 22\epsilon + 13\epsilon - |r'-s_x| \;=\; r
- 22\epsilon + 13\epsilon - r' + s_x\\
&\leq& r - 22\epsilon + 13\epsilon - r + 8\epsilon +
22\epsilon + 8\epsilon \;=\; 29\epsilon.
\eeq
This, and $|r' - r|\leq 8\epsilon$, yield
\beq
d_\infty(\Lambda(\eta x, r),\; \Lambda(x', r')) &\leq & d(x', \eta x) + |r -
r'| \;\leq\; 37\epsilon\\
\Rightarrow\quad d_\infty(\Lambda(\eta x, r),\; \kappa'\Lambda(x,r)) &\leq&
43\epsilon.
\eeq

\begin{figure}[t]
\centerline{\epsfig{file=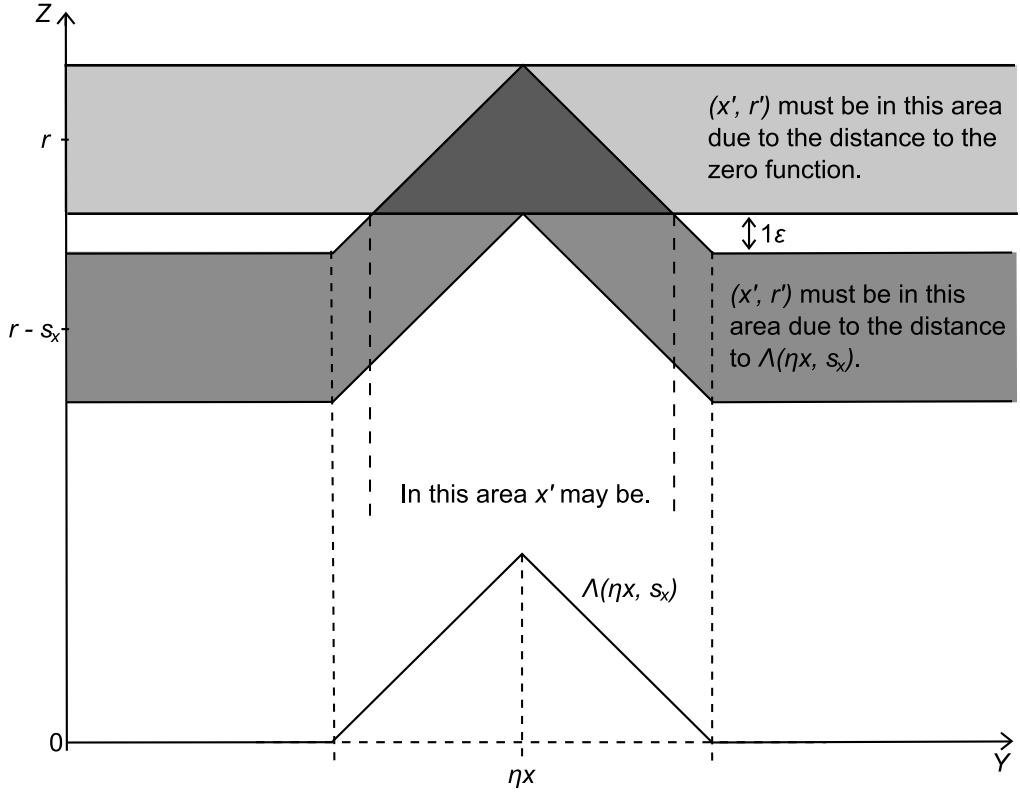}}
\caption{The function $\Lambda(x', r')$ in the proof of Lemma \ref{LEM_mt1} is
already determined up to nearness by its distance to two other functions: the
zero function and $\Lambda(\eta x, s_x)$. This shows: A $\Lambda$-function
$\Lambda(y,s)$ is not only mapped near another $\Lambda$-function
$\Lambda(y',s')$, but $y'$ only depends on $y$ and $s'$ only depends on
$s$ (modulo some multiples of $\epsilon$).} 
\label{FIG_main_theorem_2}
\end{figure}

{\bf Case 2:} $\epsilon > 0$, $r\in [0,38\epsilon)$. Obviously,
\beq
d_\infty(\Lambda(\eta x, r),\; \kappa'\Lambda(x, r)) \;\leq&&
d_\infty(\Lambda(\eta x, r), \Lambda(\eta x, s_x))\\ &+& d_\infty(\Lambda(\eta
x, s_x), \kappa'\Lambda(x,22\epsilon)) \\
&+& d_\infty(\kappa'\Lambda(x, 22\epsilon), \kappa'\Lambda(x,r))\\
\leq&& |r-s_x| + 6\epsilon + \epsilon + |r-22\epsilon|
\eeq
As $r\in [0,38\epsilon)$ and $s_x\in [14\epsilon,
30\epsilon]$ (see above), we
receive $|r-s_x|\leq 30\epsilon$ and $|r-22\epsilon|
\leq 22\epsilon$. This 
adds up to $59\epsilon$.

{\bf Case 3:} $\epsilon=0$, $r\in[0,\infty)$. As of Lemma
\ref{LEM_lambda_to_lambda}, for all $x\in X$ we can choose $\eta(x)$
such that $\kappa'\Lambda(x, 1) = \Lambda(\eta x, s_x)$ for some
$s_x\in Z\setminus \{\infty\}$. From Proposition \ref{PRO_ml-iso_zero}
we see $\kappa'(0)=0$, hence $s_x=1$. Now, let $r\in [0, \infty)$ be
arbitrary. Let $x'\in Y$, $r'\in Z$ such that $\kappa'\Lambda(x,r)
= \Lambda(x', r')$. Clearly, from the distance to 0 we again have
$r'=r$. From
\beq
d_\infty(\Lambda(x, 1), \Lambda(x, r)) &=& |r - 1|
\eeq
we conclude
\beq
d_\infty(\Lambda(\eta x, 1), \Lambda(x', r)) &=& |r - 1|.
\eeq
Due to Proposition \ref{PRO_lambda_distances} this can happen iff
(a) $|r-1|=1\geq r$ or (b) $|r-1|=r\geq 1$ or (c) $d(x',\eta x) =
0$. Case (c) prooves our statement, case (b) can't happen: $|r-1|=r$
iff $r=\frac{1}{2}$, which contradicts $r \geq 1$. So,
assume case (a). Then $r=2$, which contradicts $r \leq 1$, or
$r=0$. But the case $r=0$ is trivial, as we already saw from
Proposition \ref{PRO_ml-iso_zero} that
\beq
\kappa'\Lambda(x,0) \; = \; 0 \; = \; \Lambda(\eta x, 0).
\eeq

{\bf Case 4:} $r=\infty$. We know $\Lambda(x,\infty)=\bigsup_{s\in
[38\epsilon,\infty)} \Lambda(x, s)$. Using our result for finite $r$,
we conclude
\beq
d_\infty\left(\kappa'\bigsup_{s\in[38\epsilon,\infty)} \Lambda(x,s),\;
\bigsup_{s\in [38\epsilon,\infty)} \Lambda(\eta x, s)\right) &\leq&
1\epsilon + 43\epsilon.
\eeq
Apply $\bigsup_{s\in[38\epsilon,\infty)} \Lambda(\eta x, s) = \Lambda(\eta
x, \infty)$ to see that $\kappa'\Lambda(x,\infty)$ is
$44\epsilon$-near $\Lambda(\eta x, \infty)$.
\end{Proof}

\begin{Lemma}\label{LEM_mt2}
$\eta: X\rightarrow Y$ as defined in the proof of Lemma \ref{LEM_mt1} is a
$88\epsilon$-isometry.
\end{Lemma}
\begin{Proof}
 From Corollary \ref{COR_lambda_metric} follows
\beq
d(\eta x, \eta y) &=& \lim_{r\rightarrow \infty,\, r\neq \infty}
d_\infty\big(\Lambda(\eta x, r), \Lambda(\eta y, r)\big).
\eeq
Applying Lemma \ref{LEM_mt1} for large enough $r$ yields:
\beq
\big|d_\infty\big(\Lambda(\eta x, r), \Lambda(\eta y, r)\big) -
d_\infty\big(\kappa'\Lambda(x,r), \kappa'\Lambda(y,r)\big)\big| &\leq&
2\cdot 43\epsilon
\eeq
And of course:
\beq
\big|d_\infty\big(\kappa'\Lambda(x,r), \kappa'\Lambda(y,r)\big) -
d_\infty\big(\Lambda(x,r), \Lambda(y,r)\big)\big| &\leq& \epsilon
\eeq
Hence
\beq
\big|d(\eta x, \eta y) - d(x,y)\big| &\leq& 87\epsilon,
\eeq
i.e. $\eta$ is a rough isometric embedding. Just as $\eta$ was constructed
from $\kappa'$, we construct $\eta'$ from $\kappa$. It remains to show that
$\eta\circ\eta'$ and $\eta'\circ\eta$ are near identities. Again, we make use of
Corollary \ref{COR_lambda_metric}:
\beq
\big|d(\eta\eta' x, x) - \lim_{r\rightarrow \infty, \,r\neq\infty}
d_\infty\big(\Lambda(\eta\eta'x,r), \Lambda(x,r)\big)\big| &=& 0\\
\Rightarrow \big|d(\eta\eta' x, x) - \lim
d_\infty\big(\kappa'\kappa\Lambda(x,r),\Lambda(x,r)\big)\big|&\leq&
2\cdot 43\epsilon + \epsilon\\ 
\Rightarrow d(\eta\eta' x, x) &\leq & 88\epsilon
\eeq
Same for $\eta'\circ\eta$.
\end{Proof}

\begin{Theorem}[$=$ Th. \ref{THEP_ind_isom}]\label{THE_ind_isom}
Let $X, Y$ be complete (possibly infinite) metric spaces and $\epsilon \geq 0$.
For each $\epsilon$-ml-isomorphism $\kappa: \Lip(Y)\rightarrow \Lip(X)$ there
is a $88\epsilon$-isometry $\eta: X\rightarrow Y$, such that $\kappa$
is $61\epsilon$-near $\bar{\eta}: f\mapsto
\overline{f\circ \eta}$. 
\end{Theorem}
\begin{Proof}
Construct $\eta$ as in Lemma \ref{LEM_mt1}. It's a
$88\epsilon$-isometry due to Lemma \ref{LEM_mt2}. It remains to show
that $\kappa$ is near $\bar{\eta}$: Let $f\in \Lip Y$ be arbitrary. Represent
$f$ via $\Lambda$-functions as in Proposition \ref{PRO_lambda_building}.
Obviously, 
\beq
d_\infty\left(\kappa\bigsup_{y\in Y}\Lambda(y, fy), \bigsup_{y\in
Y}\Lambda(\eta' y, fy)\right) &\leq& 1\epsilon +
59\epsilon
\eeq
due to Lemma \ref{LEM_mt1}. Apply Proposition \ref{PRO_lambda_exchange}.
\end{Proof}

\section{Scaling limits}

\begin{Definition}
Let the {\em rough distance} $d_R(X,Y)$ between two (possibly infinite) metric
spaces $X$ and $Y$ be the infimum over all $\epsilon\geq 0$ such that $X$ and
$Y$ are $\epsilon$-isometric, or $\infty$ if there are none. If $d_R(X,Y)=0$,
the spaces $X$ and $X'$ will be called {\em pseudo-isometric}.
\end{Definition}

The rough distance fulfills triangle-inequality, as concatenation of
an $\epsilon$- and a $\delta$-isometry is an $\epsilon+\delta$-isometry. It is
closely related to the Gromov-Hausdorff-Distance for compact spaces, but may
differ in a variable between $1/2$ and $2$ (i.e., they are
Lipschitz-equivalent, see e.g. \cite{Gromov_Riemannian}, Proposition 3.5). 

Pseudo-isometry is a little bit less than isometry. However, they are
equivalent if only compact spaces are compared (e.g. \cite{Petersen},
\cite{Gromov_Riemannian}), or if we deal with simple graphs, due to their
integer metric. A nice article about scaling limits, Gromov-Hausdorff distances
and quasi-isometries in the case of graphs and Cayley graphs is
\cite{Requardt_Continuum}.

\begin{Definition}
Let $\M$ be the non-small groupoid of all pseudo-isometry-classes of metric
spaces with $\epsilon$-isometries as morphisms. $d_R$ is a (possibly infinite)
metric on $\M$ in a natural way.
\end{Definition}

Each of the components of $\M$ can be endowed with a metric
and topology, with the only drawback of being proper classes. This
``topology'' allows us to define the convergence of metric spaces to another
metric space, up to pseudo-isometry. $\M$ is complete in this ``topology''
(cf. \cite{Petersen}, Proposition 6, the proof works in non-compact and
non-separable cases as well).

\begin{Definition}\label{DEF_scaling_limit}
Let $\ell > 0$, and $s_\ell : \M \rightarrow \M$ given by
\beq
s_\ell[(X,d)] := [(X, \ell\cdot d)]
\eeq
which scales each metric space in $\M$ by the factor $\ell$ ($\ell\cdot
\infty := \infty$). This operation clearly is compatible with
pseudo-isometry. Let $[X]$ be a space in $\M$. If the limit 
\beq
s [X] := \lim_{\ell\rightarrow 0} s_\ell [X]
\eeq
exists for all sequences $\ell\rightarrow 0$, then $s[X]$ (resp. all members of
$s[X]$) is called the {\em (strong) scaling limit} of $[X]$.
\end{Definition}

We now want to apply Theorem \ref{THE_ind_ml-iso}.

\begin{Corollary}
Let $X, Y$ be some (possibly infinite) metric spaces, such that $Y$ is a strong
scaling limit of $X$ ($Y$ is unique up to pseudo-isometry). Then there is a
strong scaling limit of $\Lip X$, and it is pseudo-isometric to $\Lip Y$.
(``The scaling limit of the Lipschitz space is the Lipschitz space of the
scaling limit.'')
\end{Corollary}
\begin{Proof}
As $d_R(Y,s_\ell X) \rightarrow 0$ for $\ell \rightarrow 0$, there are
$\epsilon_\ell$-isometries $\eta_\ell:s_\ell X\rightarrow Y$ with 
$\epsilon_\ell\rightarrow 0$.
These induce $4\epsilon_\ell$-ml-isomorphisms
$\bar{\eta_\ell}: \Lip Y\rightarrow \Lip s_\ell X$, which are in
particular $4\epsilon_\ell$-isometries. Hence, $d_R(\Lip Y, \Lip s_\ell X)
\rightarrow 0$. Proper rescaling of the associated Lipschitz functions further
shows $s_\ell \Lip X$ is naturally isometric to $\Lip s_\ell X$, hence $s_\ell
\Lip X \rightarrow \Lip Y$ up to pseudo-isometry.
\end{Proof}

Note that we can restrict to a set of $\M$ when calculating a scaling
limit. Thus, we can make use of Banach's fixed point theorem if $d_R$
restricts to a true metric on this set.

\section{Perspectives}

\subsection{Generalizations}

There are several obvious ways to generalize the two main theorems: Changing the
target space $Z$ or the metric on $\Lip X$ would break the main points of the
proof, however single ideas might survive. The use of other types of functions
is a similarly difficult question:

\begin{Example}
Take $X=\{0\} \subset Y=\{0,1\} \subset \R$ and $\eta$ the inclusion, $\epsilon
= 1$. The metric spaces of $Z$-valued continuous functions $C(X)$ and $C(Y)$
with sup-norm are isomorphic to $Z$ and $Z^2$ respectively, which are not
roughly isometric.
\end{Example}

Another point is the inclusion of quasi-isometries. Although many ideas still
work in the context of quasi-isometries, a function's Lipschitz constant is
distorted in the process of Lemma \ref{LEM_lipschitzisation}. Hence there
happens to be a ``mixing'' of the Lipschitz function spaces $\Lip_K X$, which
creates deep problems and at the same time great potential: If we find a
workable solution to this problem, a new class of function spaces for groups
would emerge, ``quasi-Lipschitz functions'', so to speak.

Another very promising approach is to explore the rough isometries of
Haj\l asz-Sobolev spaces (\cite{Heinonen_Analysis}, chapter 5). These are
subsets of $L^p$ function spaces, with a norm similar to the Sobolev norm. This
norm contains a version of derivative which might compensate the obstruction we
encounter with functions of arbitrary Lipschitz constant, at least for
$p=\infty$.

\subsection{Category Interpretation}

Let $X,Y$ be (possibly infinite) metric spaces. We define
\beq
d_{ml}(\Lip X, \Lip Y) \;:=\; \inf \{\epsilon\geq 0: \exists \kappa: \Lip Y
\rightarrow \Lip X\n{ $\epsilon$-ml-isom.}\}\\
\n{ and }\quad \Lip \M \;:=\; \{\Lip X: X\n{ (poss.inf.) compl. metric
  space}\}\big/ (d_{ml}=0)
\eeq
$\Lip \M$ is well-defined as $d_R(X,Y) = 0$ iff $d_{ml}(\Lip X, \Lip Y) =
0$ and because each pseudo-isometry-class contains a complete metric space. $\Lip
\M$ is a non-small groupoid with ml-isomorphisms as morphisms. In these terms,
the mapping $\bar{\cdot}: \eta\mapsto \bar{\eta}$ is a Lipschitz equivalence
between the metric categories $\M$ and $\Lip \M$, and a contravariant functor
up to nearness of rough isometries.

\subsection{Further Remarks}

The proofs we presented here not only make use of the lattice structure of $\Lip
X$, but of a metric on it as well. In this sense, the comparison with Kaplanskys
Theorem \ref{THE_kaplansky} is inconsistent. Indeed, already a simple scaling
argument shows that we can't fully dispense with a structure besides the lattice
to reconstruct all rough isometries. Thus, how much of the coarse geometry is
really encoded in the lattice alone, and what else do we need to reconstruct
rough or quasi-isometries? E.g., does the addition of the ``Lipschitzized
scaling''
\beq
\alpha_\ell:\; \Lip X \;\rightarrow\; \Lip X, \quad f\;\mapsto\; \bigsup_{x\in
X}\Lambda(x, \ell\cdot fx)
\eeq
for $\ell \geq 0$ as a structural component already suffice?

Finally, note the similarity of 
Definition \ref{DEF_ml-iso} and the definition
of Ulam's approximate group homomorphisms in \cite{Ulam_Problems}, section VI.1;
see \cite{Hyers_Rassias} for a survey on this topic. Indeed, we can state the
question of stability of ml-homomorphisms and this directly corresponds to the
rigidity of rough isometries through our main theorems. 

$\quad$

{\bf Acknowledgements.} We want to thank the ``Graduiertenkolleg Gruppen und
Geometrie'' for supporting our research. Thanks go to Prof. Andreas Thom, Prof.
Thomas Schick, Johannes H\"artel and Dr. Manfred Requardt for interesting
discussions and hints on this subject and to Prof. Themistocles Rassias for
pointing us at Hyers' and Ulam's works.

Georg-August-Universit\"at G\"ottingen, Germany\newline
eMail \verb|lochmann@uni-math.gwdg.de|

\end{document}